\documentclass{amsart}

\def\C{{\mathcal C}}

\def\bZ{{\mathbb Z}}

\def\L{{\mathcal L}}
\def\M{{\mathcal M}}

\def\th{{\theta}}

\def\iso{{\, \cong\, }}
\def\e{{\xi}}
\def\r{{r}}

\def\<{\langle}
\def\>{\rangle}
\def\_u{{\mathfrak u}}

\def\U{\mathcal U}
\def\u{{u}}
\def\v{{v}}
\def\i{i}

\usepackage{amsmath}
\usepackage{amssymb}
\usepackage{amsfonts}

\usepackage{xy}
\input xy
\xyoption{all}

\newtheorem{theorem}{Theorem}
\newtheorem{prop}[theorem]{Proposition}
\newtheorem{lemma}[theorem]{Lemma}
\newtheorem{remark}[theorem]{Remark}
\theoremstyle{definition}
\newtheorem{definition}[theorem]{Definition}
\theoremstyle{remark}

\title{GENUS 2 FIELDS WITH DEGREE 3 ELLIPTIC SUBFIELDS}

\author{T. Shaska}

\address{Department of Mathematics,  University of California at  Irvine,  Irvine, CA  92697.}
\email{tshaska@math.uci.edu}

\begin{document}

\begin{abstract}
In this paper we study genus 2 function fields $K$ with degree 3
elliptic subfields. We show that the number of $Aut(K)$-classes of
such subfields of fixed  $K$ is 0,1,2 or 4. Also we compute an equation
for the locus of such $K$ in the moduli space of genus 2 curves.
\end{abstract}

\maketitle
\date{}

\section{Introduction}
%&&&&&&&&&&&&&&&&&&&&&&&&&&&&&&&&&&&&&&&&&&&&&&&&&&&&&&&&&&&&&&&&

We study genus two curves $\C$  whose function fields  have a degree 3 subfield of
genus 1.   Such subfields we call elliptic subfields.
 Such curves $\C$  have already occurred in the work of
 Hermite, Goursat, Burkhardt, Brioschi, and Bolza, see Krazer \cite{Krazer} (p.
 479).
More generally, degree $n$ elliptic subfields of genus 2 fields
have been studied by Frey \cite{Fr}, Frey and Kani \cite{FK}, Kuhn
\cite{Ku}, Shaska \cite{S-J}, Shaska and Voelklein
 \cite{Sh-V}.   In the degree 3 case,
 explicit equations can be used to answer questions which remain
 open in the general case. This was the theme of the author's  PhD
 thesis (see \cite{thesis})  which also covered the case of degree  2
(which is  subsumed  in \cite{Sh-V}).

Equation \eqref{eq_F1_F2} gives a normal form for pairs
$(K,E)$ where $E$ is a degree 3 elliptic subfield of $K$. This
normal form depends on two parameters $a,b\in k$. Isomorphism
classes of pairs $(K,E)$ are parameterized by parameters $u,v$ where
$u:=ab$ and $v:=b^3$.
This yields an expression for the  invariants $(i_1, i_2,
i_3)$ of $K$ (see Igusa \cite{Ig}) in terms of $u$ and $v$. Our central object of study
is the map $$\th: (u,v) \to (i_1, i_2, i_3)$$ This map is shown to have  degree
2, which means that a general genus 2 field having elliptic
subfields of degree 3 has exactly two such subfields. 
There are exactly four genus 2 fields which have four $Aut(K)$-classes and 
a 1-dimensional family of such fields with one $Aut(K)$-class of such subfields.

Although the map $\th$ is given by explicit equations (see
\eqref{eq_i}),  it is a non-trivial computational task to show that it
has degree 2 and to compute the action of its Galois group on $u$
and $v$.  The direct approach exceeds available computer power, so
we introduce auxiliary parameters $r_1$, $r_2$ that parameterize
pairs of cubic polynomials. They arise from the fact that the
subfield $E$ induces a particular sextic defining $K$ which splits
naturally as a product of 2 cubic polynomials. We show that $\th$
factorizes as $$(u,v)\to (r_1, r_2)\to (i_1, i_2, i_3)$$ where the
latter map is birational, and the former has degree 2.  Thus
$r_1$, $r_2$ yield a birational parameterization of the locus
$\L_3$ of genus 2 fields having a degree 3 elliptic subfield.
  That $\L_3$ is a rational variety follows
also from the general theory of "diagonal modular surfaces", see Kani
\cite{Kani}. We also compute the equation in $i_1, i_2, i_3$ that
defines $\L_3$ as a sublocus of the moduli space $\M_2$ of
genus 2 curves.
  We further
find relations between the j-invariants of the degree 3 elliptic
subfields of $K$ and classify all $K$ in $\L_3$ with extra
automorphisms.

All the computations were done using Maple \cite{Maple}.

\bigskip

{\bf Acknowledgment:} I would like to express my sincere gratitude to my PhD advisor 
Prof. H. Voelklein  for all the time and effort spent in guiding me towards my 
dissertation (from which this paper originated).

\bigskip

\medskip
%&&&&&&&&&&&&&&&&&&&&&&&&&&&&&&&&&&&&&&&&&&&&&&&&&&&&&&&&&&&&&&&&&&&&&&&&&&&&
\section{Genus Two Fields With Degree 3 Elliptic Subfields}

Let $k$ be an algebraically closed field of characteristic 0. All function fields 
will be over $k$.
Let $K$ be a genus 2 function field and $E$ a degree 3 elliptic
subfield of $K$.

\medskip

{\it The associated extension $k(X)/k(U)$.} It is well-known that
$K$ has exactly one genus zero subfield of degree 2, which we
denote by $k(X)$. The generator of $Gal(K/k(X))$ is the
hyperelliptic involution of $K$. It fixes each elliptic subfield
of $K$, see Tamme \cite{Tamme}. Hence the field $E\cap k(X)$,
denote it by $k(U)$, is a subfield of $E$ of degree 2.

$$\begin{matrix}
K & \buildrel{2}\over{-} & k(X) \\
3\,  | &  & | \,  3 \\
E & \buildrel{2}\over{-} & k(U)
\end{matrix}$$

\medskip

{\it Ramification of $K/E$.}  Either $K/E$ is ramified at exactly 2
places of $K$, of ramification index 2, or at one place of $K$, of
ramification index 3. It follows immediately from the Riemann-Hurwitz formula. 
The former (resp., latter) case we call the
non-degenerate (resp., degenerate) case, as in  \cite{Fr}.

\begin{definition}
A {\bf non-degenerate pair} (resp., {\bf degenerate pair}) is a pair $(K,E)$ such
that $K$ is a genus 2 field with a degree 3 elliptic subfield $E$ where
the extension $K/E$ is ramified in two (resp., one) places.  Two such
pairs $(K,E)$ and $(K' ,E')$ are called isomorphic if there is a
$k$-isomorphism $K\to K'$ mapping $E\to E'$.
\end{definition}

{\it Ramification of $k(X)/k(U)$.}  In the non-degenerate
(resp., degenerate ) case $k(X)/k(U)$ is ramified at exactly 4
(resp., 3) places of $k(X)$ each of ramification index 2 (resp., one of
index 3 and the other two of index 2), see \cite{Fr}.

\medskip

{\it Invariants of $K$.}  We denote by $J_2, J_4, J_6, J_{10}$ the
classical invariants of $K$, for their definitions see \cite{Ig}
or \cite{Sh-V}.  These are homogeneous polynomials (of degree
indicated by subscript) in the coefficients of a sectix $f(X,Z)$
defining $K$
$$Y^2=f(X,Z)=a_6X^6+ a_5X^5Z + \dots + a_1XZ^5+a_0$$
and they are a complete set of $SL_2(k)$-invariants (acting by
coordinate change).
The  absolute invariants
\begin{small}
\begin{equation}
i_1:=144 \frac {J_4} {J_2^2}, \quad i_2:=- 1728 \frac {J_2J_4-3J_6} {J_2^3}, \quad
i_3 :=486 \frac {J_{10}} {J_2^5}
\end{equation}
\end{small}
are even $GL_2(k)$-invariants.
Two genus 2 curves with $J_2\neq 0$ are isomorphic if and only if they have the same
absolute invariants.

\bigskip

{\bf Main Theorem:}
%\label{mainthm}
{\it
Let $K$ be a genus 2 field and $e_3(K)$ the number of $Aut(K/k)$-classes
 of elliptic subfields of $K$ of degree 3.  Then;

i)  $e_3(K) =0, 1, 2$, or  $4$

ii)    $e_3(K) \geq 1$ if and only if
 the classical invariants of $K$ satisfy  the irreducible
equation   $F(J_2, J_4, J_6, J_{10})=0$ displayed  in  Appendix A.
}

\medskip

There are exactly two genus 2 curves (up to isomorphism) with $e_3(K)=4$, see 4.2.  
The case $e_3(K)=1$ (resp., 2) occurs for a 1-dimensional (resp., 2-dimensional) 
 family of genus 2 curves, see section 4.

\subsection{The non-degenerate  case}

Let $(K,E)$ be a non-degenerate pair and $k(X)$  and $k(U)$ their
associated genus 0 subfields.
Both $k(X)/k(U)$ and $E/k(U)$ are ramified at 4 places of $k(U)$,
three of which are in common.  Take the common places to be $U=q_1,
q_2, q_3$. Also, take $U=0$ (resp., $U=\infty$) the place ramified in
$k(X)$ but not in $E$ (resp., in $E$ but not $k(X)$). Take $X=0$
(resp., $X=\infty$) the place over $U=0$ of ramification index 2
(resp., 1).
In the following figure bullets (resp., circles) represent
 places of ramification index 2 (resp., 1).

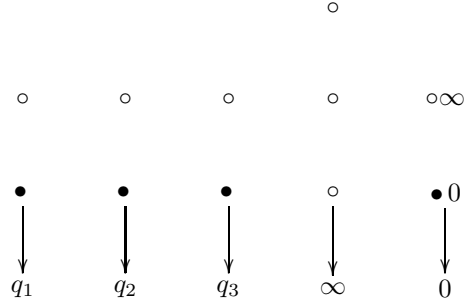
\begin{figure}[!hbt]
$$
\xymatrix{& & & \circ   \\
\circ    & \circ    & \circ   & \circ   & \circ \infty  \\
\bullet \, \ar[d]  &\bullet \,   \ar[d] &\bullet \,   \ar[d] &  \circ
  \ar[d]  &\bullet \,  0  \ar[d] \\
q_1 &q_2 &q_3 &  \infty  & 0      \\
}
$$
\caption{Ramification of $k(X)/k(U)$}
\label{fig2}
\end{figure}

Thus,
$U=\phi(X)=l\frac {X^2} {X^3+a X^2+b X+c}$,
where $l, a,b,c\in k$ such that $X^3+a
X^2+b X+c\, $ has no multiple roots and $l,c \neq 0$.
We can make $l=1$ and $c=1$ by
replacing  $U$ and $X$ by scalar multiples.
Then,
\begin{equation}\label{phi1}
U=\phi(X)=\frac {X^2} {X^3+a X^2+b X+1}
\end{equation}

\begin{remark}\label{rem_1}
These conditions determine $X$ up to multiplication by a third root of
unity $\e_3\in k$. Replacing $X$ by $\e_3X$ replaces $(a,b)\,  $  by $\, \,  (\e_3^2 a,
\e_3 b)$. Invariants of this transformation are:

\begin{equation}\label{u_v}
\begin{split}
 u\, = \, ab\\
v\, =\, b^3\\
\end{split}
\end{equation}

\end{remark}

\medskip

\noindent The derivative of $\phi(X)$  is
$$\phi^\prime (X)=-\frac {X(X^3-bX-2) } {(X^3+a X^2+b X+1)^2 } $$
Taking the resultant with respect to $z$  of  $z^2-bz-2$ and $\frac
{\phi(x)-\phi(z)} {x-z}$ we get
\begin{equation}\label{Z_eq}
(x^3-bx-2)(4x^3+b^2x^2+2bx+1)
\end{equation}
The roots of this polynomial correspond to the 6 places of $k(X)$ over
the places $U=q_1, q_2, q_3$. Thus, the 6 places of $k(X)$ ramified in
$K$ correspond to the roots of the polynomial on the right side of
\begin{equation}\label{eq_F1_F2}
Y^2=(X^3+a X^2+b X+1 )\, (4X^3+b^2X^2+2bX+1 )
\end{equation}
which gives an equation of $K$.
The discriminant of the
sextic is nonzero, hence
\begin{equation}
\Delta:=(4a^3+27-18ab-a^2b^2+4b^3)^2(b^3-27) \neq 0
\end{equation}

One checks that the element
\begin{equation}\label{V}
V= Y\,  \frac {X^3-bX-2} {F(X)^2}
\end{equation}
satisfies
$$V^2=(U-q_1)(U-q_2)(U-q_3)= U^3 +2 \frac {ab^2-6a^2+9b} {R} \, U^2 +
\frac {12a-b^2 } {R} \, U - \frac 4 {R} $$
where $R:=\frac {\Delta} {(b^3-27)}$.
Thus, $E=k(U,V)$.

\begin{lemma}\label{lemma1}  For $(a,b)\in k^2$ with $\Delta\neq 0$, equation
\eqref{eq_F1_F2} defines a  genus 2 field  $K_{a,b}=k(X,Y)$. It has a
non-degenerate
degree 3 elliptic subfield $E_{a,b}=k(U,V)$, where $U$ and $V$ are given in
\eqref{phi1} and \eqref{V}.
Two such  pairs $(K_{a,b}, E_{a,b})$ and  $(K_{a', b'},
E_{a', b'})$ are isomorphic if and only if $u=u'$ and $v=v'$
(where $u,v$ and $u',v'$ are associated with $a,b$ and
$a', b'$, respectively, by \eqref{u_v}).
\end{lemma}

\proof
The first two statements follow  by reversing the above
arguments. If $u=u'$, $v=v'$ then $a'=\e_3^i a$, $b'=\e_3^i b$ for some
$i$. Then clearly the two associated non-degenerate pairs are isomorphic. The
converse follows from Remark \ref{rem_1}.

\qed
%-----------------------------------------------------------------------

\begin{prop}
The $(u,v)\in k^2$ with $\Delta\neq 0$  bijectively parameterize the
isomorphism classes of non-degenerate pairs   $(K,E)$ (via the
parameterization defined in Lemma 3).
\end{prop}

\proof The proof follows from Lemma 3.

\qed

%\begin{remark}
\noindent From  the normal form of $K$ in \eqref{eq_F1_F2}
we can compute the
classical invariants $J_2, J_4, J_6, J_{10}$  in
terms of $u,v$. These expressions satisfy equation
\eqref{eq_L3_2} which proves one implication of claim ii) of the
theorem for the non-degenerate pairs. In section 4 we explain how
equation \eqref{eq_L3_2} was found.
%\end{remark}

%&&&&&&&&&&&&&&&&&&&&&&&&&&&&&&&&&&&&&&&&&&&&&&&&&&&&&&&&&&
\subsection{The degenerate case}
%&&&&&&&&&&&&&&&&&&&&&&&&&&&&&&&&&&&&&&&&&&&
Let  $K/E$ be a degenerate pair. Then
$k(X)/k(U)$ is ramified at  exactly three places of $k(U)$ which are
also ramified in $E$. Let $U=0$ be the place ramified in $E$ but not
in $k(X)$. Take $U=\infty$ be the place of $k(U)$ that is totally
ramified in $k(X)$ and $X=\infty$ the place over it.

\begin{figure}[!hbt]
$$\xymatrix{& & & \circ   \\
   & \circ   & \circ  & \circ       \\
\buildrel{\infty}\over\odot  \ar[d]  &\bullet \,   \ar[d] &\bullet \,   \ar[d] &  \circ
  \ar[d]  & \\
\infty & & &  0        \\
}
$$
\caption{Ramification of $k(X)/k(U)$, degenerate case.}
\label{fig_deg}
\end{figure}
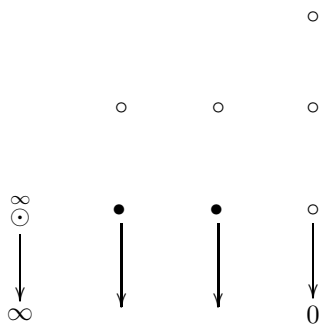

Then, $U=l(X^3+aX^2+bX+c)$ for $a, b, c, l\in k$ such that
$X^3+aX^2+bX+c$  has no
multiple roots and  $l\neq 0$. Replacing $X$ by $X+\frac a 3$ and $U$
by a scalar multiple we get
$U=X^3+bX+c$.
Note that $b\neq 0$ (otherwise  two of the places of $k(U)$ ramified
in $k(X)$ coalesce).  Thus, by replacing $X$ by a scalar multiple we
can also get $b=1$. Then,
\begin{equation}
U=X^3+X+c
\end{equation}
As in Lemma \ref{lemma1} we get
 $K=k(X,Y)$ such that
\begin{equation}\label{deg_case}
Y^2=(3X^2+4)(X^3+X+c)
\end{equation}
where $c^2\neq - \frac 4 {27}$.
Then, $V=Y(3X^2+1)$
satisfies the equation
$$V^2=U(27U^2-54c U+4+27c^2)$$
Thus, $E=k(U,V)$.

\begin{remark} In this case, those  $w:=c^2$ with $w\neq - \frac 4
{27} $  bijectively parameterize the isomorphism classes of degenerate pairs
$K/E$. The proof is analogous to the non-degenerate case.
\end{remark}

\noindent Since
\begin{equation}
\begin{split}
J_2& =774\\
J_4&=36(268+837w)\\
J_6 &= 36(76760 +290574w -729w^2)\\
J_{10} &=432(27w+4)\\
\end{split}
\end{equation}
\noindent the map $$k\setminus \{- \frac 4 {27} \} \to k^3$$
$$w \to (i_1, i_2, i_3)$$
is injective. Thus, a genus 2 field $K$ has at most one $Aut(K)$-orbit
of subfields $E$ such that $(K,E)$ is  a degenerate pair.
 One can check that $J_2, J_4, J_6, J_{10}$ satisfy equation
\eqref{eq_L3_2}. This completes the proof that the condition in part (ii) of 
Main Theorem is necessary.

>From (10) we compute that the locus of those $K$ being part of a 
degenerate pair $(K,E)$ is given by

\begin{small}
\begin{equation}\label{eq_deg_case}
\begin{split}
 & 3418801i_1^2-2550732480i_1+4023336960+611249816i_2=0\\
&161446939560824832i_3+23630752512i_1^2-6321363049i_1^3-29445783552i_1\\
& +12230590464=0
\end{split}
\end{equation}
\end{small}

\noindent (Note $J_2=774\neq 0$ on this locus, so $i_1, i_2, i_3$ are everywhere defined). 
For later use we note the following: If $(K,E_0)$ is a non-degenerate
pair with parameters $u, v$  then there exists a degenerate pair $(K,E_1)$ (for same $K$) if
and only if
\begin{equation}\label{sp_uv}
2v-9u+27=0
\end{equation}
This is obtained by expressing $i_1,i_2,i_3$ in \eqref{eq_deg_case} in
terms of $u,v$.

%%%%%%%%%%%%%%%%%%%%%%%%%%%%%%%%%%%%%%%%%%%%%%%%%%%%%%

\section{Function Field of $\L_3$}
%&&&&&&&&&&&&&&&&&&&&&&&&&&&&&&&&&&&&&&&&&&&&&&&&&&&&&&&&&&&&&&&&&&
The absolute invariants  $i_1, i_2, i_3$  in terms of $u,v$ are
\begin{small}
\begin{equation}\label{eq_i}
\begin{split}
i_1 = &\frac {144}  {v(-405+252u+4u^2-54v-12uv+3v^2)^2}  (1188u^3-8424uv+u^4v-24u^4\\
& +14580v -66u^3v+138uv^2+297u^2v+945v^2-36v^3+9u^2v^2)  \\
i_2 =&- \frac {864} {v^2(-405+252u+4u^2-54v-12uv+3v^2)^3}
(-81v^3u^4+2u^6v^2+234u^5v^2\\
& +3162402uv^2-21384v^3u+26676v^4-473121v^3-72u^6v -5832v^4u+14850v^3u^2\\
& -72v^3u^3+324v^4u^2-650268u^3v-5940u^3v^2-3346110v^2+432u^6-1350u^4v^2\\
& +136080u^4v -7020u^5v-307638u^2v^2   \\
i_3 =& -243 \frac {(v-27)(4u^3-u^2v-18uv+4v^2+27v)^3}{v^3(-405+252u+4u^2-5
4v-12uv+3v^2)^5} \\
\end{split}
\end{equation}
\end{small}
We will sometimes view $u,v$ as parameters from $k$ and sometimes as
the corresponding coordinate functions on $k^2$. From the context it
will be clear which point of view we are taking.

\begin{lemma}
$[k(i_1,i_2,i_3):k(u,v)]\leq 2$
\end{lemma}

\proof
>From the resultants of  equations in \eqref{eq_i} we determine
that $[k(\v):k(\i_1,\i_2)]=16$, $[k(\v):k(\i_2,\i_3)]=40$,  and
$[k(\v):k(\i_1,\i_3)]=26$.
 We  can show that $\u\in k(\i_1,
\i_2, \i_3, \v)$, the expression is large and we don't display it
(the interested reader can check \cite{thesis}).
Since $[k(\u,\v):k(\i_1,\i_2, \i_3)]$ must be a common divisor of 16,
 26, and 40, then the claim follows.

\qed

We need to show that $[k(\u,\v):k(\i_1,\i_2, \i_3)] = 2$. Since we
know that $[k(u,v):k(i_1,i_3)]=26$, this is equivalent to showing that
$[k(i_1,i_2,i_3):k(i_1,i_3)]=13$. This follows once we have equation
\eqref{eq_L3_2} in terms of $i_1,i_2,i_3$ but at this stage we cannot
derive it: Eliminating $u,v$ from  equations
 \eqref{eq_i}   exceeds available computer power.
 We use additional invariants $\r_1, \r_2$
  to overcome this problem.

%&&&&&&&&&&&&&&&&&&&&&&&&&&&&&&&&&&&&&&&&&&&&&&&&&
\subsection{Invariants of Two Cubics}
%&&&&&&&&&&&&&&&&&&&&&&&&&&&&&&&&&&&&&&&&&&&&&&&&&&&&&&&&&&&&&&&&&&&&&&&&&&&
We define the following  invariants of two cubic polynomials.
For $F(X)=a_3X^3+a_2 X^2 + a_1 X+a_0$ and $G(X)=b_3
X^3+b_2X^2+b_1X+b_0$
define
$$H(F,G) :=a_3 b_0 - \frac  1 3 a_2 b_1 + \frac 1 3 a_1 b_2 -a_0 b_3$$
We denote by $R(F,G)$ the resultant of $F$ and $G$ and  by $D(F)$ the
discriminant of $F$.
Also,
$$r_1(F,G) =\frac  {H(F,G)^3} {R(F,G)}, \quad r_2(F,G)=\frac
{H(F,G)^4} {D(F)\, D(G)}$$

\begin{remark}
In \cite{Vishi} it is shown that $r_1, r_2$, and $ r_3=\frac {H(F,G)^2} {J_2(F\,G)}$ form a
complete system of invariants for unordered pairs of cubics.
\end{remark}

\noindent For 
$F(X)=X^3+aX^2+bX+1      $ and 
$G(X)=4X^3+b^2X^2+2bX+1$  and $u,v$ as in \eqref{u_v}
 we have
\begin{equation}
\begin{split}\label{eq_r}
r_1(F,G) &= 27\frac {v(v-9-2u)^3} {4v^2-18uv+27v-u^2v+4u^3}\\
r_2 (F,G)& = -1296 \frac  {v(v-9-2u)^4}   {(v-27)(4v^2-18uv+27v-u^2v+4u^3)}\\
\end{split}
\end{equation}

Taking the resultants from the above equations we get the following
equations for $\u$ and $\v$ over $k(\r_1, \r_2)$:

\begin{small}
\begin{equation}
\begin{split} \label{eq_u}
65536\r_1\r_2^3 \,\u^2+(42467328\r_2^4+21233664\r_2^4\r_1+480\r_2\r_1^4+2\r_1^5+41472\r_2^2\r_1^3\\
+1548288\r_2^3\r_1^2-294912\r_2^3\r_1)\u-382205952\r_2^4+238878720\r_2^4\r_1-2654208\r_2^3\r_1\\
+13934592\r_2^3\r_1^2+285696\r_2^2\r_1^3+ 2400\r_2\r_1^4+7\r_1^5 & =0\\
\end{split}
\end{equation}
\end{small}
\begin{small}
\begin{equation}
\begin{split} \label{eq_v}
16384\v^2\r_2^3+(221184\r_2^3\r_1+\r_1^4+11520\r_2^2\r_1^2-442368\r_2^3+192\r_2\r_1^3)\v\\
-5971968\r_2^3\r_1-864\r_2\r_1^3-124416\r_2^2\r_1^2-2\r_1^4 & =0 \\
\end{split}
\end{equation}
\end{small}

%&&&&&&&&&&&&&&&&&&&&&&&&&&&&&&&&&&&&&&&&&&&&&&&&&&&&&&&&&&&&&&&&&&&&&&&&&&&&

Roots of  equation  \eqref{eq_u}  (resp., equation \eqref{eq_v})   are $\u$ and
$\beta(\u)$ (resp., $v$ and $ \beta(v)) $ where,

\begin{small}
\begin{equation}\label{eq_nu}
\begin{split}
\beta(\u) & =\frac {(\v-3\u)(324\u^2+15\u^2\v-378\u\v-4\u\v^2+243\v+72\v^2)}
{(\v-27)(4\u^3+27\v-18\u\v-\u^2\v+4\v^2)}\\
\beta(\v) & =- \frac {4(\v-3\u)^3}{4\u^3+27\v-18\u\v-\u^2\v+4\v^2} \\
\end{split}
\end{equation}
\end{small}

It follows that either $[k(u,v):k(r_1,r_2)]=2$ or $k(u,v)/k(r_1,r_2)$
is  Galois with Klein 4-group as Galois group. In the latter case, the
map
\begin{equation}
\begin{split}
\tau: k(\u,\v) & \to k(\u, \v)\\
u    &  \to \u \\
v      & \to \beta(\v) \\
\end{split}
\end{equation}
\noindent is an involutory automorphism of $k(\u,\v)$. But by plugging in
equation  \eqref{eq_nu} we see that it is not involutory. Thus, $[k(\u,\v):k(\r_1,\r_2)]=2$ and
$Gal_{k(\u,\v)/k(\r_1, \r_2)}=\<\beta\>$.

%&&&&&&&&&&&&&&&&&&&&&&&&&&&&&&&&&&&&&&&&&&&&&&&&&&&&&&&&&&&&&&&&&&
\begin{lemma} The fields   $k(\i_1, \i_2, \i_3)$  and $k(\r_1, \r_2)$
coincide,  hence $\, \, [k(u,v):k(i_1,i_2,i_3)]=2$. Moreover;
\begin{small}
\begin{equation}\label{i1_i2_i3}
\begin{split}
\i_1 &= \frac 9 4 \frac {(13824\r_1^3\r_2^2+442368\r_1^2\r_2^3+5308416\r_1\r_2^4+192\r_1^4\r_2+\r_1^5+786432\r_1\r_2^3+9437184\r_2^4)}{\r_1(-1152\r_2^2+96\r_2\r_1+\r_1^2)^2} \\
\i_2 &=\frac {27} {8\r_1^2(-1152\r_2^2+96\r_2\r_1+\r_1^2)^3 } (+79626240\r_1^4\r_2^4-4076863488\r_1^2\r_2^5+34560\r_1^6\r_2^2\\
& +12230590464\r_1^2\r_2^6+32614907904\r_1\r_2^6+14495514624\r_2^6 +288\r_1^7\r_2+2211840\r_1^5\r_2^3\\
& +\r_1^8-212336640\r_1^3\r_2^4+1528823808\r_1^3\r_2^5 -2359296\r_1^4\r_2^3)\\
\i_3  &=-521838526464 \frac {\r_2^9} {\r_1^2(-1152\r_2^2+96\r_2\r_1+\r_1^2)^5}\\
\end{split}
\end{equation}
\end{small}
\end{lemma}

\proof
Equations \eqref{i1_i2_i3}   are verified by expressing all variables in terms
of $u,v$ (These equations were found by taking suitable resultants). It
follows that $k(i_1,i_2,i_3) \subset k(r_1, r_2)$.  Equality follows
since $[k(u,v):k(\r_1,\r_2)]=2$ and
 $[k(\u,\v):k(\i_1,\i_2,\i_3)]\leq 2$ (see Lemma 6).

\qed

\begin{remark}
To find equation \eqref{eq_L3_2} we  eliminate $\r_1$ and
$\r_2$ from  equations \eqref{i1_i2_i3}.
 The resulting equation  \eqref{eq_L3_2} has degree 8, 13, and 20 in 
$\i_1, \i_2, \i_3$ respectively.
\end{remark}

\section{Proof of Main Theorem}
%&&&&&&&&&&&&&&&&&&&&&&&&&&&&&&&&&&&&&&&&&&&&&&&&&&&&&&&&&&
%\subsection{The general case}
The map
$$\th: (u,v) \to (i_1, i_2, i_3)$$
given by \eqref{eq_i}  has degree 2, by previous section and it is
defined when $J_2\neq 0$. For now we assume that $J_2\neq 0$ (The
case $J_2=0$ is treated in section 4.2). Denote the minors of the
Jacobian matrix of   $\th$ by $M_1(u,v), M_2(u,v), M_3(u,v)$. The
solutions of

\begin{equation}
\begin{split}
\left\{
\aligned
  M_1(u,v)= 0 \\
  M_2(u,v)=0\\
  M_3(u,v)=0\\
\endaligned
\right.
\end{split}
\end{equation}
consist of the (non-singular) curve
\begin{equation}
\begin{split}\label{n_3_iso1}
8v^3+27v^2-54uv^2-u^2v^2+108u^2v+4u^3v-108u^3=0\\
\end{split}
\end{equation}
and  7 isolated solutions which we display in table \eqref{tab1}
together with the corresponding  values $(i_1, i_2, i_3)$ and
properties of the corresponding genus 2 field $K$.
%&&&&&&&&&&&&&&&&&&&&&&&&&&&&&&&&&&&&&&&&&&&&&&&&&&&&&&
\begin{small}
\begin{table}[!ht]
\renewcommand\arraystretch{1.5}
\noindent\[
\begin{array}{|c|c|c|c|c|}
\hline
(u,v) & (i_1, i_2, i_3) & Aut(K) & e_3(K) \\
\hline
(-\frac 7 2, 2) & J_{10}=0, \quad \text{no associated} &  &  \\
 & \text{genus 2 field K}  &  & \\
\hline
(-\frac {775} 8, \frac {125} {96}), & & & \\
 (\frac {25} 2, \frac {250} {9})&
- \frac {8019} {20}, -\frac {1240029} {200}, \frac {531441} {100000}     & D_4 & 2 \\
\hline
(27- \frac {77} 2 \sqrt{-1}, 23+ \frac {77} 9 \sqrt{-1}), &  & & \\
(27+ \frac {77} 2 \sqrt{-1}, 23- \frac {77} 9 \sqrt{-1}) & (\frac
{729} {2116}, \frac {1240029}  {97336}, \frac {531441}
{13181630464}  & D_4 & 2   \\
\hline
(-15+ \frac {35} 8 \sqrt{5}, \frac {25} 2 + \frac {35} 6 \sqrt{5}), &  & & \\
(-15- \frac {35} 8 \sqrt{5}, \frac {25} 2 - \frac {35} 6 \sqrt{5})&
81, - \frac {5103} {25}, -\frac {729} {12500}  &  D_6 & 2  \\
\hline
\end{array}
\]
\caption{Exceptional points where $det(Jac(\th))=0$}
\label{tab1}
\end{table}
\end{small}

Let $\bar u, \bar v, \bar i_1, \bar i_2, \bar i_3$ denote the restrictions of the corresponding 
functions to the curve \eqref{n_3_iso1}. Eliminating $u$ from \eqref{n_3_iso1} and the 
defining equation for $i_1$ in (13) we get a relation between $v$ and $i_1$ which 
shows $[k(\bar v, \bar i_1):k(\bar i_1)]=8$; also we get $\bar u$ expressed as a rational function of 
$\bar v$ and $\bar i_1$. 
Hence $k(\bar v, \bar i_1)=k(\bar u, \bar v)$ and so $[k(\bar u, \bar v):k(\bar i_1)]=8$.  
Similarly we get a relation between $\bar v$ and $\bar i_2$. Eliminating $\bar v$ gives 
\begin{equation}\label{curve}
G(\bar i_1, \bar i_2)=0
\end{equation}
an irreducible equation between $\bar i_1$ and $\bar i_2$ of degree 8 (resp., 12) 
in $\bar i_1$ (resp., $\bar i_2$).  It is fully displayed as equation \eqref{eq_J} in Appendix B. 
Thus, $[k(\bar i_1, \bar i_2):k(\bar i_1)]=8$.  
Since also $[k(\bar u, \bar v):k(\bar i_1)]=8$, it follows that $k(\bar u, \bar v)=k(\bar i_1, \bar i_2)$. 
Hence,
$e_3(K) = 1$ for any $K$ such that the associated $u$ and $v$
satisfy equation \eqref{n_3_iso1} 
except possibly those $(u,v)$ that lie over  singular points of the curve 
\eqref{curve}.  We check that all these singular points 
 have multiplicity 2, hence there are at most two 
$(u,v)$ points over each of them. 

\begin{remark}
If $K$ has one subfield of degree 3, then it has at least two distinct such subfields.
This follows 
from the explicit construction in \cite{thesis}, Lemma 5.4. (It is actually true for degree $n$ 
elliptic subfields, for any $n$, see \cite{Sh-V}).  Since $e_3(K)=1$ for a generic curve $K$ of the locus 
\eqref{curve}, such $K$ has $|Aut(K)| > 2$. Actually $Aut(K)=V_4$. We verified this by using 
equation (17) from 
\cite{Sh-V} which gives a necessary and sufficient condition for a genus 2 field $K$ to have $V_4 \leq Aut(K)$
(in terms of the classical invariants of $K$). Taking the resultant of this equation and the one in Appendix A
 yields an expression that has \eqref{curve} (the projectivized version)  as a factor. From the conditions 
given in \cite{Sh-V} for $K$ to have $Aut(K)$ strictly bigger than $V_4$, one checks that actually $Aut(K)=V_4$ for all 
$K$ in locus \eqref{curve}. 
\end{remark}

\subsection{The general case} 
%**************************************
By previous section $\th$ is generically a covering of
degree 2. So there exists a Zariski open subset $\U$ of $k^2$ with
the following properties: Firstly, $\th$ is defined everywhere on
$\U$ and is a covering of degree 2 from $\U$ to $\th(\U)$.
Further, if $\_u \in \U$ then all $\_u^\prime \in k^2$ with $\th$
defined at $\_u^\prime$ and $\th(\_u^\prime)=\th(\_u)$ also lie in
$\U$. Suppose $\underline i \in k^3$ such that $|
\th^{-1}(\underline i)|
>2$ and $det(Jac(\th))$ does not vanish at any point of $\th^{-1}(\underline i) $.
Then by implicit function theorem, there is an open ball $B$ around
each element of $\th^{-1}(\underline i) $ such that each point in
$\th(B)$ has $> 2$ inverse images under $\th$. But $B$ has to
intersect the Zariski open set $\U$. This is a contradiction. Thus, if
${\underline i} \in k^3$ and $|\th^{-1}(\underline i)| > 2$, then
$det (Jac(\th)) =0$ at some point of $\th^{-1}(\underline i) $ and so
$\underline i$ satisfies (22). 
It follows that if $K$ is a genus 2 field with $J_2\neq 0$ and
$G \neq 0$ then $e_3(K)= 2$.

Note that if there is a degenerate pair $(K,E_1)$ and a
non-degenerate pair $(K,E_2)$ (for the same $K$) then for the
$(u,v)$-invariants of the latter we have \eqref{sp_uv}, hence
$\beta(v)=27$ and so the $\beta$-image of $(u,v)$ doesn't correspond
to a non-degenerate pair.

We can now complete the proof of part (ii) of Main Theorem.  We noted  before that the condition in
(ii) is necessary. This condition (see Appendix A) defines an irreducible sublocus of $\mathcal M_2$ because 
the corresponding equation is irreducible (by Maple). 
The locus $e_3(K) \geq 1$ lies in this sublocus and contain a 
dense subset of it (image of $\theta$). By Lange \cite{Lange}, the locus $e_3(K)\geq 1$ is Zariski closed, 
hence equals the locus defined by the equation in Appendix A. 
\medskip
%************************

%&&&&&&&&&&&&&&&&&&&&&&&&&&&&&&&&&&&&&&&&&&&&&&&&&
\subsection{Exceptional Cases for $J_2=0$.}
%&&&&&&&&&&&&&&&&&&&&&&&&&&&&&&&&&&&&&&&&&&&&&&&&&&&&&&&&&&&&&&&&&&&&&&&&&&&&&
We consider the case when $J_2=0$, where

\begin{equation}
\begin{split}
J_2 & =-2(3v^2+4u^2-12uv+252u-54v-405)\\
 &= \frac 2 3 (3v-6u+2u\sqrt{6}-27-18\sqrt{6})(3v-27+18\sqrt{6}-6u-2u\sqrt{6})
\end{split}
\end{equation}

Genus 2 fields $K$ with $J_2=J_4=0$ and $J_6\neq 0$ (resp.,  $J_2=J_6=0$ and
  $J_4\neq 0$)   are classified up to isomorphism by  the invariant $\frac {J_6^5} {J_{10}^3}$
(resp.,  $\frac {J_4^5} {J_{10}^2}$),  see Igusa \cite{Ig}. In the first  (resp., second) case 
 there are exactly two (resp., four) such $K$ with $e_3(K) \geq 1$ and
they all have  $e_3(K)=2$.

The invariants
\begin{equation}\label{eq_a_1_a_2}
a_1:= \frac {J_4 \cdot J_6} {J_{10}}, \quad a_2:=\frac {J_6 \cdot J_{10}} {J_4^4}
\end{equation}

\noindent determine genus two fields with $J_2=0$, $J_4\neq 0$, and $J_6\neq 0$  up to isomorphism,
see Igusa \cite{Ig}. A field $K$ with such invariants has $e_3(K) \geq 1$ if and only if

\begin{tiny}
\begin{equation}\label{J2_curve_1}
\begin{split}
46656a_1^5a_2^3+7558272a_1^4a_2^3+1259712\sqrt{6}a_1^4a_2^3-15552a_1^4a_2^2
-12427478784a_1^3a_2^2+4917635712\sqrt{6}a_1^3a_2^2\\
+1728a_1^3a_2-656217531654480a_1^2a_2^2
+267571209034080\sqrt{6}a_1^2a_2^2-1844125056a_1^2a_2+743525568\sqrt{6}a_1^2a_2\\
-64a_1^2-6334497449472117312a_1a_2^2+2585860435265558832\sqrt{6}a_1a_2^2+230833239838992a_2a_1-\\
94237227087840a_2a_1\sqrt{6}-601244429975805030777a_2^2+245429539257764380572a_2^2\sqrt{6} &=0
\end{split}
\end{equation}
\end{tiny}

or

\begin{tiny}
\begin{equation}\label{J2_curve_2}
\begin{split}
230833239838992a_2a_1+94237227087840a_2a_1\sqrt{6}-12427478784a_1^3a_2^2-4917635712\sqrt{6}a_1^3a_2^2\\
+1728a_1^3a_2-15552a_1^4a_2^2-656217531654480a_1^2a_2^2-267571209034080\sqrt{6}a_1^2a_2^2+46656a_1^5a_2^3-64a_1^2\\
-601244429975805030777a_2^2-245429539257764380572a_2^2\sqrt{6}-1844125056a_1^2a_2-743525568\sqrt{6}a_1^2a_2\\
-6334497449472117312a_1a_2^2-2585860435265558832\sqrt{6}a_1a_2^2+7558272a_1^4a_2^3-1259712\sqrt{6}a_1^4a_2^3 &=0
\end{split}
\end{equation}
\end{tiny}

As one can check by discussing the map 
$$\vartheta: (u,v)\to (a_1,a_2)$$
these fields $K$ generically have $e_3(K) = 2$. It remains to check   the singular
points of the curves \eqref{J2_curve_1} and \eqref{J2_curve_2}.
The intersection of the two  components is empty.
 Each component has four singular points two of which have $e_3(K)=2$ and the other two have $e_3(K)=4$.  
Thus there are exactly four genus 2 fields with $e_3(K) = 4$.  Two of them  satisfy \eqref{J2_curve_1},
 we call them $P_1$ and $P_2$. The other two are obtained from $P_1, P_2$ by 
applying the automorphism $\sqrt{6} \to - \sqrt{6}$.

\begin{table}[!ht]
\renewcommand\arraystretch{1.5}
\noindent\[
\begin{array}{|c|c|}
\hline
P_1 & a_1=-\frac {77169}{8}+\frac {30759}{8}\sqrt{6}, \quad 
a_2=\frac {13783592}{23149125}+\frac {5629912}{23149125}\sqrt{6}   \\
\hline
P_2 &   a_1=-\frac {650835}{8}+\frac {268785}{8}\sqrt{6}, \quad 
a_2=-\frac {2984}{20002028625}-\frac {144872}{60006085875}\sqrt{6}  \\
\hline
\end{array}
\]
\caption{Singular points of \eqref{J2_curve_1} with $e_3(K)=4$}
\label{tab_J2_0}
\end{table}

The four  values of $(u,v)$ corresponding to $P_1$ (resp., $P_2$)
are obtained as follows: $u$ is determined by $v$ via

$$u=\frac 1 2 (3+\sqrt{6})(v-9-6\sqrt{6})$$

and  $v$ is one of the 4 roots of 

\begin{small}
\begin{equation}\label{eq_1}
\begin{split}
9v^4+(-693-27\sqrt{6})v^3+(15141+1749\sqrt{6})v^2+(-66414-33396\sqrt{6})v\\
-1368\sqrt{6}+3388 &=0
\end{split}
\end{equation}
\end{small}

(resp.)

\begin{small}
\begin{equation}\label{eq_2}
\begin{split}
3v^4+(-153+33\sqrt{6})v^3+(3129-1219\sqrt{6})v^2+(-10854-11556\sqrt{6}) v\\
+176904\sqrt{6}+450036 &=0
\end{split}
\end{equation}
\end{small}

This completes the proof of the Main Theorem.

\qed

%&&&&&&&&&&&&&&&&&&&&&&&&&&&&&&&&&&&&&&&&&&&&&&&&&&&&&&&&&&&&&&&&
\section{j-invariants}
%&&&&&&&&&&&&&&&&&&&&&&&&&&&&&&&&&&&&&&&&&&&&&&&&&&&&&&&

A genus 2 field $K$ corresponding to a generic point of locus $e_3(K) \geq 1$ has exactly 2 
elliptic subfields $E_1$ and $E_2$ of degree 3. We can take $E_1=E$ from   lemma \ref{lemma1}. 
Its j-invariant is 
\begin{small}
\begin{equation}\label{j_1}
\begin{split}
j_1 &= 16v \frac {(vu^2+216u^2-126vu-972u+12v^2+405v)^3}
{ (v-27)^3(4v^2+27v+4u^3-18vu-vu^2)^2}\\
\end{split}
\end{equation}
\end{small}

The automorphism $\beta \in Gal_{k(u,v)/k(\r_1, \r_2)}$ permutes
 $E_1$ and $E_2$, therefore switches $j_1$ and
$j_2$. Then $j_2=\beta(j_1)$ is
\begin{small}
\begin{equation}\label{j_2}
\begin{split}
j_2 &= -256 \frac {(u^2-3v)^3}{v(4v^2+27v+4u^3-18vu-vu^2)} \\
\end{split}
\end{equation}
\end{small}
For the  equation of $E_2$ using a different approach see \cite{thesis}, chapter 5.

\subsection{Isomorphic Elliptic Subfields}
%&&&&&&&&&&&&&&&&&&&&&&&&&&&&&&&&&&&&&&&&&&&&&&&&&&&&&&&&&&&&&&&&&&&&&&&&&&&&
Suppose that $E_1\iso E_2$. Then, $j_1=j_2$ implies that
\begin{small}
\begin{equation}
\begin{split}
8v^3+27v^2-54uv^2-u^2v^2+108u^2v+4u^3v-108u^3=0\\
\end{split}
\end{equation}
\end{small}
or
\begin{small}
\begin{equation}
\begin{split}\label{n_3_iso2}
& 324v^4u^2-5832v^4u+37908v^4-314928v^3u-81v^3u^4+255879v^3+30618v^3u^2\\
& -864v^3u^3-6377292uv^2 +8503056v^2-324u^5v^2+2125764u^2v^2-215784u^3v^2\\
& +14580u^4v^2+16u^6v^2+78732u^3v+8748u^5v -864u^6v-157464u^4v+11664u^6 =0\\
\end{split}
\end{equation}
\end{small}
\begin{remark}
The former equation is the condition that $det ( Jac(\th))=0$ see
equation (21).
\end{remark}

\begin{remark}
If $e_3(K)\geq 1$ then the automorphism group
$Aut(K/k)$  is one of the following: $\bZ_2, V_4$, $ D_8$, or $D_{12}$.
$\bZ_2$ is the generic case, $V_4$ is a 1-dimensional family.
There are exactly 6 genus 2 fields with automorphism
group $D_8$ (resp., $D_{12}$).
\end{remark}

For the proof see \cite{thesis}. 

\bigskip

\begin{center}
{\bf Appendix A:}
\end{center}
%&&&&&&&&&&&&&&&&&&&&&&&&&&&&&&&&&&&&&&&&&&&&&&&&&&&&&&&&&&&&&&&&&&&&&&&&&&&&&

\bigskip

\noindent Here is the equation which defines the locus  $\L_3$ of fields $K$ with $e_3(K) \geq 1$. 
\begin{equation}
\label{eq_L3_2}
C_8 J_{10}^8 +C_7 J_{10}^7+ \dots +  C_1J_{10} +C_0=0
\end{equation}
where  $C_0, \dots , C_8$ are

\medskip

\begin{tiny}
$C_8  = 2^4\cdot 3^{31\cdot} 5^5\cdot {19}^{10}\cdot {29}^5$
\end{tiny}

\begin{tiny}
$C_7 = {2}^4\cdot 3^{27} \cdot {19}^5(-194894640029511J_2^5-55588661819356000J_4^2J_2- 12239149540657725J_2^3J_4\\
 \quad  \quad      +
 223103526505680000J_4J_6+40811702108053500J_2^2J_6)
$
\end{tiny}

\medskip

\begin{tiny}
$C_6 = 2^2\cdot 3^{21} (-35802284468757765858432J_4^5-1756270399106587730391J_4^2J_2^6-28638991859523006654J_4J_2^8\\
 -84091225203760159441286J_4^3J_2^4+400895959391006953561032J_4^4J_2^2-61773685738999443J_2^10\\
 -3673201396072259603756160J_4^3J_2J_6+7879491755218264984387200J_4^2J_6^2+15251447355608658629952J_2^5J_4J_6\\
 +1179903008384844066250272J_2^3J_4^2J_6-5566672398589809889658760J_2^2J_4J_6^2+112024289372554183680J_2^7J_6\\
 -32116769409722716182888J_2^4J_6^2+8512171877754962249155200J_2J_6^3)
$
\end{tiny}

\medskip

\begin{tiny}
$
C_5 =2^2 \cdot 3^{19}(-12630004382695462653J_4^4J_2^7+320839252764287362560J_4^7J_2-1876069272397136886448J_4^6J_2^3\\
 +606742866220456356580J_4^5J_2^5-124173485719052715J_4^3J_2^9+22241034512101438944000J_4^5J_2^2J_6\\
 -88546736703938826304512J_4^4J_2J_6^2-10712078926420753449984J_2^4J_4^4J_6+68904635323303664511264J_2^3J_4^3J_6^2\\
 -192353895694677164016384J_2^2J_4^2J_6^3+197449733923926905783808J_2J_4J_6^4-1916173047371645223936J_4^6J_6\\
 +116211018774997425051648J_4^3J_6^3-132143597697786172416J_6^5+1361403542457288J_2^10J_4J_6\\
 -5005765118740492656J_2^7J_4J_6^2+232819061639483430720J_2^6J_4^3J_6-1576319894694585178452J_2^5J_4^2J_6^2\\
 +4655239459208764553088J_2^4J_4J_6^3-226900590409548J_2^{11}J_4^2-2042105313685932J_2^9J_6^2+6261632755967872800J_2^6J_6^3\\
 -5065734796478275576176J_2^3J_6^4+1352318109350828796J_2^8J_4^2J_6)
$
\end{tiny}

\medskip

\begin{tiny}
$
C_4 = 3^{15}(1417825317153277312J_4^9J_2^2+2391308818408811717J_4^6J_2^8+718590303030600J_2^{10}J_4^5\\
   -638760745337170544640J_4J_6^6+440759275303802880J_4^{10}-8118717280771686540192J_4^5J_2^4J_6^2\\
 +42668434906863398019072J_4^2J_2J_6^5-57054664814020640574336J_4^3J_2^2J_6^4+30546774740158581676032J_4^4J_2^3J_6^3\\
 +8601814215123831275904J_4^6J_2^2J_6^2-1449562700682195916800J_4^7J_2^3J_6+1067928354124249303104J_4^6J_2^5J_6\\
-10443896263024316301312J_2^3J_4J_6^5+7247970315150439028112J_2^4J_4^2J_6^4-21769241176751736619008J_4^5J_2J_6^3\\
 -2640201919999154595648J_2^5J_4^3J_6^3-55893562424445261312J_2^7J_4^5J_6+531409635241191119304J_2^6J_4^4J_6^2\\
 -1012614205133520J_2^8J_6^4-2454855015326199552J_2^5J_6^5-12501409939920J_2^{12}J_4^4-675076136755680J_2^{10}J_4^2J_6^2\\
 -33390518666828400J_2^9J_4^4J_6+3188363568027498432J_2^6J_4J_6^4-1569001498547402304J_2^7J_4^2J_6^3\\
 -275375222428239820800J_4^7J_6^2+19809849095518050330624J_4^4J_6^4+6179516061983740183680J_2^2J_6^6\\
 +150016919279040J_2^{11}J_4^3J_6+87799481406335621136J_4^8J_2^4+1350152273511360J_2^9J_4J_6^3\\
 -55496611186132800648J_4^7J_2^6+39911809855842557952J_4^8J_2J_6+353362680242481096J_2^8J_4^3J_6^2   )
$
\end{tiny}

\medskip

\begin{tiny}
$
C_3  = 2^4\cdot 3^{12} (-19225816442103600J_4^{10}J_2^5+6433952690394144J_2^4J_6^7-2917203075615J_2^{11}J_4^7\\
+62951605613640J_2^{10}J_4^6J_6 +7900854051362368J_4^{11}J_2^3-873165547551982J_2^9J_4^8+13234982161044480J_4^{12}J_2\\
+4077902864550187008J_4^5J_6^5 +7506792545698293J_4^9J_2^7-55019014994202624J_4^{11}J_6-3415519987075510272J_4^2J_6^7\\
-932605137272623104J_4^8J_6^3 -6607177263254292480J_2J_6^8-1394785406520J_2^7J_6^6-1913285880J_2^{13}J_4^6\\
-258293593800J_2^{11}J_4^4J_6^2 -1976299597616301504J_4^8J_2^3J_6^2+1337598192058041744J_4^7J_2^5J_6^2\\
-2324642344200J_2^9J_4^2J_6^4 +2789570813040J_2^8J_4J_6^5-243015467955111198J_2^7J_4^6J_6^2+22136761801348668J_2^8J_4^7J_6\\
-155463896437263612J_4^8J_2^6J_6 +16101033796183004352J_4^5J_2^3J_6^4+16367298631796450304J_4^3J_2J_6^6+8433152J_4^7J_2^2J_6^3\\
-8254965178021469184J_4^6J_2J_6^4 +34439145840J_2^{12}J_4^5J_6-576988130682378J_2^9J_4^5J_6^2+2912934238489260J_2^8J_4^4J_6^3\\
-8749875412454175J_2^7J_4^3J_6^4 +15637511592200340J_2^6J_4^2J_6^5-127105068829245696J_4^{10}J_2^2J_6+614908581517421568J_4^9J_2J_6^2\\
-23374419431360207616J_4^4J_2^2J_6^5 +1508868948605946984J_2^6J_4^5J_6^3-5795040294470623824J_2^5J_4^4J_6^4\\
+14094983896511630112J_2^4J_4^3J_6^5+1033174375200J_2^{10}J_4^3J_6^3+314069798204069472J_4^9J_2^4J_6-61501104J_4^6J_2^4J_6^3\\
 -21194163080222025024J_2^3J_4^2J_6^6+18002402119176332544J_2^2J_4J_6^7-15392091937240080J_2^5J_4J_6^6)
$
\end{tiny}

\medskip

\begin{tiny}
$
C_2 =2^5\cdot 3^8(-159732958548480J_4^{13}J_2J_6-27945192968593920J_2J_4J_6^9+238596124150086J_4^{10}J_2^7J_6+3224288J_2^3
J_6^9\\
 -36311136215244J_2^9J_4^9J_6-996173640J_2^{12}J_4^6J_6^2+5977041840J_2^{11}J_4^5J_6^3-22413906900J_2^{10}J_4^4J_6^4-375014140095
08J_4^{13}J_2^4\\
 +14210312049697149J_2^6J_4^6J_6^4+86354918885580768J_2^4J_4^4J_6^6-111444977082978432J_2^3J_4^3J_6^7+27908893977856J_4^{14}J_2^2\\
 +83768141083825152J_2^2J_4^2J_6^8-42942980968765488J_2^5J_4^5J_6^5-10736445647473J_4^{11}J_2^8+61746352553318400J_4^4J_2J_6^7\\
 +410958880454688J_4^{12}J_2^3J_6-14059252057660416J_4^3J_6^8-1643659809866496J_4^{11}J_2^2J_6^2+441832778741790J_2^8J_4^8J_6^2\\
 -3128599551108636J_2^7J_4^7J_6^3+2815950495430656J_4^{10}J_2J_6^3+12291244885171152J_4^8J_2^5J_6^3+1579225145J_2^{12}J_4^9\\
 -2286353789913249J_4^9J_2^6J_6^2-40300476525629352J_4^7J_2^4J_6^4+81707043798929088J_4^6J_2^3J_6^5-25936092270J_2^6J_6^8\\
 -98257765274489088J_4^5J_2^2J_6^6+3318887207480832J_4^{10}J_2^4J_6^2-478511899451856J_4^{11}J_2^5J_6+28476287051677J_4^12J_2^6\\
 -12153253649302656J_4^9J_2^3J_6^3+24757975700165376J_4^8J_2^2J_6^4-26570902457981952J_4^7J_2J_6^5-6755065089024J_4^{15}\\
 +10883911680J_6^{10}+53793376560J_2^9J_4^3J_6^5-80690064840J_2^8J_4^2J_6^6+69162912720J_2^7J_4J_6^7+94873680J_2^{13}J_4^7J_6\\
 +227109129291J_2^{10}J_4^7J_6^2-628213747356J_2^9J_4^6J_6^3-1389130574661J_2^8J_4^5J_6^4+16465793988870J_2^7J_4^4J_6^5\\
 -56794191944715J_2^6J_4^3J_6^6+102713329135152J_2^5J_4^2J_6^7-98529746457492J_2^4J_4J_6^8-30650938650J_2^{11}J_4^8J_6\\
 -1716480768J_4^9J_6^4+11718053954519040J_4^6J_6^6+220752428322816J_4^{12}J_6^2+1322792799725J_2^{10}J_4^{10}-3953070J_2^{1
4}J_4^8)
$
\end{tiny}

\medskip

\begin{tiny}
$
C_1 =-2^8\cdot 3^5(61736960J_4^8J_2-182135808J_4^7J_6+16021872J_4^7J_2^3-211022400J_4^6J_2^2J_6-26594919J_4^6J_2^5\\
 +899159040J_4^5J_2J_6^2+330458928J_4^5J_2^4J_6-215198J_2^7J_4^5-1227405312J_4^4J_6^3-1535734368J_4^4J_2^3J_6^2\\
 +2930532J_4^4J_2^6J_6-363J_2^9J_4^4+3162070656J_4^3J_2^2J_6^3-16471998J_4^3J_2^5J_6^2+4356J_4^3J_2^8J_6-19602J_4^2J_2^7J_6^2\\
 -2433162240J_4^2J_2J_6^4+47961936J_4^2J_2^4J_6^3+39204J_4J_2^6J_6^3-72369936J_4J_2^3J_6^4+746496J_4J_6^5-29403J_2^5J_6^4\\
 +45116352J_2^2J_6^5)(J_4^3-J_2^2J_4^2+6J_2J_6J_4-9J_6^2)^3 \\
$
\end{tiny}

\vspace{-.5cm}
\begin{tiny}
$C_0=2^8(768J_4^2-416J_4J_2^2-J_2^4+1536J_2J_6)(J_4^3-J_2^2J_4^2+6J_2J_6J_4-9J_6^2)^6$
\end{tiny}

\bigskip

\bigskip

%******************************************************************************************

\begin{center}
{\bf Appendix B:}
\end{center}

\bigskip

\noindent The equation of the branch locus of   the map 
$$\theta: k^2 \setminus \{\Delta =0\} \to \L_3$$
$$(u,v) \to (i_1, i_2, i_3)$$

\begin{tiny}
\begin{equation}\label{eq_J}
\begin{split}
&\, \, 3507505273398025 i_1^{12}-4880484817793862073548480i_1^{11}+192(11302504938388489628125 i_2  \\
&+346452039237689650581192)i_1^{10}-4(89439905046278964319358976+60307157046030532997225i_2^2 \\
&-4188981066069113234648640 i_2)i_1^9+192(4924672355809568004885504+33588887753890413143515i_2^2\\
&-89290373874540245356608i_2)i_1^8-192(6360139591235383381327872+593587805135845078438632i_2^2\\
&+849397678885114696829952i_2-15537701670163340329775i_2^3)i_1^7+2(406606074742962841916289024i_2^2\\
&-27711519511099200420652160i_2^3+305526006347596356290347008+183049808606955774794075i_2^4\\
&+255320765313220782576893952i_2)i_1^6-192i_2(13654023711946692280289280i_2+60377228453350507376315 i_2^3\\
&-1464350022771442997265792i_2^2+2245754530466537929703424) i_1^5 +192i_2^2(19906166147692916436664320\\
&+524465846991778455117432i_2^2-2333360138590692031481856i_2-10694743521252049023175 i_2^3) i_1^4\\
&+20 i_2^2(2535928320048654288481984 i_2^3-5763664926963183845376000i_2-28262691566657845506883584i_2^2\\
&-99413657833844087193600000-9280304409903257402325i_2^4)i_1^3 -192 i_2^3(1807602326421361731982656i_2^2\\
&-76401752022567738591625i_2^3-2867627019656613888000000-9918952715424066995911680i_2)i_1^2\\
&+ 1728i_2^4(1262371230708245434125i_2^3-44970318919363276752280i_2^2-1752882366993587712000000\\
&+506870449554602235884544i_2)i_1+i_2^4(127385395640432909375625i_2^4-6059765135286837968377600i_2^3\\
&-717148324858373259264000000i_2+102838387547772598641352704i_2^2+1617386925920256\cdot
10^{12}) =0
\end{split}
\end{equation}
\end{tiny}

%&&&&&&&&&&&&&&&&&&&&&&&&&&&&&&&&&&&&&&&&&&&&&&&&&&&&&&&&&&&&&&&&&&&

\end{document}